\documentclass[12pt,a4paper,twoside]{article}
\usepackage{amsfonts}
\date{}
\newcommand{\A}[1]{\vspace{3mm}}
\newcommand{\D}{{\rm d}}

\newcommand{\R}{\mathbb{R}}

\begin{document}

\author{Rolf Schneider and Franz E. Schuster}
\title{Rotation Equivariant Minkowski Valuations}

\maketitle

\begin{center}
{\em Dedicated to Professor Peter Gruber\\on the occasion of his sixty-fifth birthday}
\end{center}

\begin{abstract} \noindent The projection body operator $\Pi$,
which associates with every convex body in Euclidean space
$\R^n$ its projection body, is a continuous valuation, it
is invariant under translations and equivariant under
rotations. It is also well known that $\Pi$ maps the set of
polytopes in $\R^n$ into itself. We show that
$\Pi$ is the only non-trivial operator with these properties. \\[2mm]
\end{abstract}
\section{Introduction and Main Results}

Let ${\cal K}^n$ denote the space of convex bodies
(non-empty, compact, convex sets) in \linebreak $n$-dimensional
Euclidean space $\R^n$ ($n\ge 2$), endowed with the
Hausdorff metric. A convex body $K \in \mathcal{K}^n$ is
determined by its support function $h(K,\cdot)$, defined on
$\mathbb{R}^n$ by $h(K,x)=\max\{ \langle x,y \rangle :y \in K\}$,
where $\langle\cdot,\cdot\rangle$ is the scalar product of
$\R^n$. The {\em projection body} $\Pi K$ of $K$ is
defined by
\[
h(\Pi K,u) = V_{n-1}(K|u^{\perp}) \qquad\mbox{for }u\in S^{n-1}.
\]
Here, $K|u^{\perp}$ denotes the image of $K$ under orthogonal
projection to the $(n-1)$-dimensional subspace orthogonal to $u$, and $S^{n-1}$ is the unit sphere of
$\R^n$.  Generally, we denote by $V_k(M)$ the $k$-dimensional volume of a $k$-dimensional
convex body $M$.

The projection body operator was already introduced by Minkowski
\cite{Min11}. In recent years it has attracted increased attention due to its
numerous applications in different areas, see \cite{bolk69,
bourg88, Gar95, grin99, Lud02,
Thom96}. Projection bodies of convex bodies are centered
convex bodies called {\em zonoids}. For their role in geometry,
we refer to the surveys \cite{SW83, GW93}. Projection
bodies of convex polytopes, called {\em zonotopes}, have appeared in
optimization, computational geometry, and other
areas, see \cite{Zieg95}.

In this paper, the emphasis is on the fact that the projection body operator
$\Pi:{\cal K}^n \to {\cal K}^n$ is a {\em Minkowski valuation},
i.e., a valuation with respect to Minkowski addition on ${\cal
K}^n$. In general, a mapping $\varphi:{\cal K}^n\to A$ into an
abelian semigroup $(A,+)$ is called a {\em valuation} if
\[ \varphi(K\cup M) + \varphi(K\cap M) = \varphi(K) + \varphi(M)\]
whenever $K,M,K\cup M\in{\cal K}^n$. Valuations on convex bodies
are a classical concept. Probably the most famous result in this
area is Hadwiger's classification of rigid motion invariant
real valued continuous valuations, see \cite{Had57,
KlaRo97} and the surveys \cite{MS83}, \cite{McM93}. In
recent years, many new results on real and body valued valuations
have been obtained, see \cite{Ales99, Ales01,
Kla95, Kla99, Lud02, Lud03, Lud05, Lud06, Lud99, Sch96}.

An immediate consequence of a result obtained by M. Ludwig in \cite[Corollary 2.2]{Lud05},
extending a previous result from \cite{Lud02},
is the following characterization of the projection body operator.

\A

\noindent{\bf Theorem.} {\em Let $\Phi:{\cal K}^n\to{\cal K}^n$
be a continuous, translation invariant valuation with the
property that, for all $K\in {\cal K}^n$ and every $\alpha\in
SL(n)$,
\[ \Phi\alpha K = \alpha^{-T}\Phi K.\]
Then there is a constant $c\ge 0$ such that $\Phi = c\Pi$.}

\A

Thus, among all continuous, translation invariant valuations
from ${\cal K}^n$ to ${\cal K}^n$, the projection body operator
is characterized, up to a factor, by its $SL(n)$
contra{\-}variance. It was also shown in \cite{Lud05} that
the assumption of continuity can be omitted when ${\cal K}^n$
as the domain of $\Phi$ is replaced by ${\cal P}^n$, the set of
convex polytopes in $\R^n$.

In the following, we will consider continuous, translation
invariant valuations $\Phi: \mathcal{K}^n \rightarrow
\mathcal{K}^n$, but we will replace the strong assumption of $SL(n)$
contravariance, which belongs to affine geometry, by the
Euclidean condition of rotation equivariance, i.e., the property
that, for all $K\in {\cal K}^n$ and every $\vartheta$ in the
rotation group $SO(n)$ of $\R^n$,
\[\Phi\vartheta K =\vartheta\Phi K.\]
The projection body operator is no longer characterized by
these properties. Simple further examples are the trivial
maps ${\cal I}$ and $-{\cal I}$ given by
\[ {\cal I}(K)= K-s(K) \quad \mbox{and}\quad (-{\cal I})(K)=-K+s(K)
\quad \qquad\mbox{for }K\in{\cal K}^n.\]
Here, $s: \mathcal{K}^n \rightarrow \mathbb{R}^n$ denotes the
Steiner point map, defined by
\begin{equation} \label{defsteinerpoint}
s(K)=n\int_{S^{n-1}} h(K,u)u \,\D u,
\end{equation}
where the integration is with respect to the rotation invariant probability
measure on the sphere. The Steiner point map is the unique vector
valued, rigid motion equivariant and continuous valuation on
$\mathcal{K}^n$, see \cite[Satz 2]{Sch72}.

A large class of non-trivial examples is provided by translation
invariant Minkowski endomorphisms. This class of operators was
introduced and investigated by the first author, see
\cite{Sch74a, Sch74}, and more recently studied by Kiderlen
\cite{Kid05}. They are precisely the continuous valuations from
${\cal K}^n$ to ${\cal K}^n$,
invariant under translations and equivariant under rotations, that
are homogeneous of degree one. Here a function $\varphi$ from
${\cal K}^n$ to $\R$ or ${\cal K}^n$ is called {\em
homogeneous of degree $j$} if $\varphi (\lambda K) =
\lambda^j\varphi (K)$ for $K\in{\cal K}^n$ and $\lambda\ge 0$.
The case of valuations homogeneous of degree $n - 1$, called
Blaschke Minkowski homomorphisms, was investigated recently by
the second author in \cite{Schu06, Schu06a}.

The main object of this paper is to find an additional assumption
which suffices to single out among this large class of valuations
the combinations of the projection body operator $\Pi$ and the
mappings $\mathcal{I}$ and $-\mathcal{I}$. This additional
assumption will be the property that polytopes are mapped to
polytopes.

\A

\noindent{\bf Theorem 1.} {\em Let $n\ge 3$. Let $\Phi:{\cal K}^n
\to {\cal K}^n$ be a continuous, translation invariant and
rotation equivariant valuation. If $\Phi$ maps polytopes to
polytopes, then
\[ \Phi = c_1\Pi + c_2 {\cal I}+c_3(-{\cal I}) \]
with constants $c_1,c_2,c_3\ge 0$.}

\A

\noindent{\bf Remarks.} The assumption that $\Phi$ maps polytopes
to polytopes is convenient to formulate, but stronger than necessary.
As the proof shows, only the following is needed. To every $j \in
\{0, \ldots, n\}$, there exists a convex body $K$ of dimension
$j$ such that $\Phi mK$ is a polytope, for $n+1$ different values of $m$.

In the plane, where the rotation group is abelian, the assertion
has to be modified. Let $\Phi:{\cal K}^2 \to {\cal K}^2$ be a
continuous, translation invariant and rotation equivariant
valuation. If the image of $\Phi$ contains some polygon with more
than one point, then there are rotations
$\vartheta_1,\dots,\vartheta_r$ of $\R^2$ and positive
numbers $\lambda_1,\dots,\lambda_r$ such that
\[ \Phi K = \lambda_1\vartheta_1 [K-s(K)]+\dots+\lambda_r\vartheta_r [K-s(K)] \]
for all $K\in{\cal K}^2$. This was proved in \cite[Satz 3]{Sch74}.

Under the additional assumption of homogeneity of degree one, the
combinations of $\mathcal{I}$ and $-\mathcal{I}$ were
characterized by the first author in \cite[Corollary 1.12]{Sch74a}.

\A

Among a subclass of the Blaschke Minkowski homomorphisms (which
includes the even ones), the projection body operator was
characterized (up to a factor) by the second author in \cite[Theorem 5.3]{Schu06}, by
the assumption that it maps some $n$-dimensional convex body to a
polytope. The wish to generalize this characterization has led us
to the following result.

\A

\noindent{\bf Theorem 2.} {\em Let $n\ge 3$. Let $\Phi:{\cal K}^n
\to {\cal K}^n$ be a continuous, translation invariant and
rotation equivariant valuation. If $\Phi$ maps bodies of
dimension $n-2$ to $\{0\}$ and maps some $n$-dimensional convex
body to a polytope, then $\Phi=c\Pi$ with some constant $c\ge 0$.}

\A

Without the assumption that $\Phi$ maps some $n$-dimensional
convex body to a polytope, every Blaschke Minkowski homomorphism
would satisfy the conditions of Theorem 2. Thus, in particular,
the multiples of the projection body operator are the only
Blaschke Minkowski homomorphisms that map some $n$-dimensional
convex body to a polytope, see Corollary 1.

\section{Results on homogeneous valuations}

In this section we collect further material on convex
bodies and some well-known results from the theory of real valued
valuations. General references are the books by
Schneider \cite{Sch93} and by Klain and Rota \cite{KlaRo97}. At the end
of this section we will prove the main tool for the proofs of
Theorems 1 and 2.

A convex body $K \in \mathcal{K}^n$ is uniquely determined by its
support function $h(K,\cdot)$, which is positively homogeneous of
degree one and sublinear. Conversely, every function with these
properties is the support function of a convex body. From the
definition of $h(K,\cdot)$ it is easily seen that $h(\vartheta
K,u)=h(K, \vartheta^{-1}u)$ for every $u \in \R^n$ and every
$\vartheta \in SO(n)$. The support function $h(K,\cdot)$ of a
convex body $K \in \mathcal{K}^n$ is piecewise linear if and only
if $K$ is a polytope.

A convex body $K \in \mathcal{K}^n$ with non-empty interior is also
determined up to translation by its surface area measure
$S_{n-1}(K,\cdot)$. For a Borel set $\omega \subseteq
S^{n-1}$, the value $S_{n-1}(K,\omega)$ is the $(n-1)$-dimensional Hausdorff measure of the set of
all boundary points of $K$ at which there exists a normal vector
of $K$ belonging to $\omega$. The relation $S_{n-1}(\lambda
K,\cdot)=\lambda^{n-1}S_{n-1}(K,\cdot)$ holds for all $K \in
\mathcal{K}^n$ and $\lambda \geq 0$. For $\vartheta \in SO(n)$,
we have $S_{n-1}(\vartheta K,\cdot)=\vartheta S_{n-1}(K,\cdot)$,
where $\vartheta S_{n-1}(K,\cdot)$ is the image measure of
$S_{n-1}(K,\cdot)$ under the rotation $\vartheta$. By Minkowski's
existence theorem, a non-negative measure $\mu$ on $S^{n-1}$ is
the surface area measure of a convex body if and only if $\mu$
has its center of mass at the origin and is not concentrated
on any great subsphere.

We collect some auxiliary results on translation invariant real
valued valuations, which will be employed repeatedly.

\A

\noindent{\bf Lemma 1} (Hadwiger \cite[p. 79]{Had57}). {\em If
$\varphi: {\cal K}^n \to \R$ is a continuous,
translation invariant valuation, homogeneous of degree $n$, then
$\varphi = cV_n$ with a constant $c$. }

\A

\noindent{\bf Lemma 2} (McMullen \cite{McM77}). {\em Every
continuous, translation invariant valuation \linebreak $\varphi:
{\cal K}^n \to \R$ has a unique representation
\[ \varphi = \varphi_0 + \dots + \varphi_n,\]
where $\varphi_j: {\cal K}^n \to \R$ is a continuous,
translation invariant valuation which is homogeneous of degree
$j$.}

\A

A valuation $\varphi$ on $\mathcal{K}^n$ is called
{\em simple} if $\varphi(K)=0$ whenever $\dim K < n$.
A function $\varphi$ from $\mathcal{K}^n$ to $\R$ or
$\mathcal{K}^n$ is called {\em even} (resp. {\em odd}\,) if
$\varphi(-K)=\varphi (K)$ (resp. $\varphi(-K)=-\varphi(K)$) for
all $K\in{\cal K}^n$. The following classification of translation invariant,
continuous and simple valuations, due to Klain \cite{Kla95} (for even
valuations) and the first author \cite{Sch96}, will be useful.

\A

\noindent{\bf Lemma 3.} {\em If $\varphi: {\cal K}^n \to {\mathbb
R}$ is a continuous, translation invariant and simple valuation,
then
\[\varphi(K)=cV_n(K)+\int_{S^{n-1}}g(u)\,\D S_{n-1}(K,u) \qquad\mbox{for }K\in \mathcal{K}^n, \]
where $c$ is a constant and $g$ is an odd, continuous real
function on $S^{n-1}$. }

\A

From Lemma 3 we deduce a result on homogeneous valuations which
are not necessarily simple. For a subspace $E$ of ${\mathbb R}^n$,
we denote by ${\cal K}(E)$ the set of convex bodies contained in $E$, and
by $SO(E)$ the subgroup of rotations in $SO(n)$ mapping $E$ into itself.
The map $\pi_E:\R^n\to E$ is the orthogonal projection.

\A

\noindent{\bf Lemma 4.} {\em Let $\varphi: {\cal K}^n \to {\mathbb
R}$ be a continuous, translation invariant valuation which is
homogeneous of degree $j$, for a given
$j\in\{0,1,\dots,n-1\}$.} \\[2mm]
(a) {\em If $\varphi$ is even and if $\varphi(K) = 0$ whenever $\dim K=j$, then $\varphi=0$.}\\[2mm]
(b) {\em If $\varphi(K) = 0$ whenever $\dim K = j + 1$, then
$\varphi=0$.}

\A

\noindent {\em Proof.} Assertion (a) was proved by Klain \cite[Corollary 3.2]{Kla99}.
In order to prove assertion (b), we can assume
$j\in\{0,1,\dots,n-2\}$. Let $E$ be a $(j+2)$-dimensional linear
subspace of $\R^n$. By the continuity of
$\varphi$, we have $\varphi(K)=0$ if $\dim K\le j+1$, thus the restriction of the
valuation $\varphi$ to ${\cal K}(E)$ is simple. It is continuous and
invariant under translations of $E$ into itself, therefore we deduce from Lemma 3, applied
in $E$, that it is a linear combination of valuations that are
homogeneous of degrees $j+2$ and $j+1$, respectively. Since
$\varphi$ is homogeneous of degree $j$, we get $\varphi=0$ on
${\cal K}(E)$. Since $E$ was an arbitrary $(j+2)$-dimensional subspace, we
have $\varphi(K)=0$ whenever $\dim K\le j+2$. Now we can repeat
the argument with a $(j+3)$-dimensional subspace, and so on, to
conclude finally that $\varphi(K)=0$ holds for all convex bodies
$K\in{\cal K}^n$. \hfill  $\Box$

\A

If $\varphi$ satisfies all assumptions of Lemma 4 (a) except the
evenness, then it follows that the valuation defined by $K\mapsto
\varphi(K)+\varphi(-K)$ is identically zero. In particular, one
can deduce that $\varphi(K)=0$ for all centrally symmetric convex
bodies $K$.

Lemma 4 leads to the following auxiliary result on valuations
taking their values in the space of convex bodies. The equation
$\Phi=\{0\}$ means that $\Phi$ maps every convex body to
the one-pointed set containing only the origin of $\R^n$.

\A

\noindent{\bf Lemma 5.} {\em Let $\Phi: {\cal K}^n \to {\cal K}^n$
be a continuous, translation invariant  and rotation equivariant valuation which is
homogeneous of degree $j$, for a given $j\in\{0,1,\dots,n-1\}$.
If $\Phi K=\{0\}$ whenever $\dim K=j$, then $\Phi =\{0\}$.}

\A

\noindent {\em Proof.} If the assumptions are satisfied, we deduce
from Lemma 4 (a) (applying it to $h(\Phi K+\Phi(-K), u)$
with $u\in\R^n$) that $\Phi$ is odd and that $\Phi K=\{0\}$ holds for all
centrally symmetric bodies $K\in{\cal K}^n$.

To extend the latter result to general convex bodies, we use an
argument employed by Klain \cite{Kla95}. Let $E$ be a
$(j+1)$-dimensional linear subspace of $\R^n$. Let
$\Delta$ be a simplex in $E$, say $\Delta={\rm conv}
\{0,v_1,\dots,v_{j+1}\}$, without loss of generality. Let
$v:=v_1+\dots+v_{j+1}$ and $\Delta':= {\rm conv}\{v,v-v_1,\dots
,v-v_{j+1}\}$. The parallelepiped, $P$, that is spanned by
$v_1,\dots,v_{j+1}$, is the union of $\Delta,\Delta'$ and a
centrally symmetric polytope $Q$, where $\dim(\Delta\cap
Q)=\dim(\Delta'\cap Q)=j$ and $\Delta\cup Q$ is convex. By assumption,
the valuation $h(\Phi(\cdot),u)$, for given
$u\in\R^n$, vanishes on convex bodies of dimension smaller
than $j+1$. Therefore, the restriction of $h(\Phi(\cdot),u)$ to
${\cal K}(E)$ is a simple valuation. It follows that
$h(\Phi\Delta,u) +h(\Phi Q,u)+ h(\Phi\Delta',u)= h(\Phi P,u)$.
Since $u\in\R^n$ was arbitrary and $\Phi Q = \{0\}$, this
yields $\Phi\Delta +\Phi\Delta'=\Phi P=\{0\}$. Since the summands
on the left-hand side are convex bodies, this is only possible if
$\Phi\Delta$ is one-pointed. Now we use a standard argument: every
polytope can be decomposed into simplices, and a simple valuation on polytopes
has an additive extension to the finite unions of polytopes.
Together with the continuity of $\Phi$, this yields that
$\Phi K$ is one-pointed for all $K\in{\cal K}(E)$, say $\Phi K=\{t^K\}$. Consider the
odd map $t: \mathcal{K}(E) \rightarrow E$ defined by $t(K):=\pi_E t^K$.
Then $t + s: {\cal K}(E) \rightarrow E$ is a continuous valuation,
equivariant under the rigid motions of $E$. By the uniqueness
property of the Steiner point map mentioned after definition
(\ref{defsteinerpoint}), this yields $t(K)=0$ for all $K \in
\mathcal{K}(E)$. Thus, $t^K$ is contained in $E^{\perp}$ for all
$K \in \mathcal{K}(E)$. If $\dim E \leq n - 2$, the vector $t^K$
is invariant under every rotation of $\mathbb{R}^n$ that
leaves $E$ pointwise fixed, thus $t^K=0$. If $\dim
E=n-1$, we have $t^K= \varphi(K)e$, where $e$ is a unit normal
vector of $E$ and $\varphi: \mathcal{K}(E) \rightarrow \R$ is an
odd, translation invariant, continuous and simple valuation. For
every $\vartheta \in SO(E)$ fixing $e$ we get
\[\varphi(\vartheta K)=h(\Phi \vartheta K,e)=h( \vartheta \Phi K,e)
=h( \vartheta \Phi K,\vartheta e)=h(\Phi K,e)=\varphi( K). \]
Thus, $\varphi$ is also invariant with respect to rotations of $E$
into itself. The valuation $\varphi$ can be represented as in Lemma 3
(applied in $E$); here $c=0$ since $\varphi$ is odd, and by the
rotation invariance, the odd function $g$ is constant and hence
vanishes. This gives $\varphi=0$ and thus $\Phi K=\{0\}$, for all $K \in \mathcal{K}(E)$.

Since $E$ was an arbitrary $(j+1)$-dimensional subspace, we have
$\Phi K=\{0\}$ for all convex bodies $K$ with $\dim K\le j+1$.
Therefore, by Lemma 4 (b), the valuation $h(\Phi(\cdot),u)$, for
given $u \in \mathbb{R}^n$, vanishes identically. Since $u$ was
arbitrary, this yields $\Phi=\{0\}$. \hfill  $\Box$

\A

In \cite{Sch74a}, the first author started an investigation of
continuous maps $\Psi: \mathcal{K}^n \rightarrow \mathcal{K}^n$,
called (Minkowski) endomorphisms, with the following properties:\\[2mm]
(a) $\Psi$ is Minkowski additive, i.e., $\Psi(K + L)=\Psi K + \Psi
L$ for $K, L \in \mathcal{K}^n.$\\[2mm]
(b) $\Psi(\vartheta K + t) = \vartheta \Psi K+t$ for $K \in
\mathcal{K}^n, \vartheta \in SO(n)$ and $t \in
\mathbb{R}^n$.\\[2mm]
Note that continuity and Minkowski additivity imply that
$\Psi$ is homogeneous of degree one. Moreover, since $(K \cup
L)+(K \cap L)=K + L$ whenever $K, L, K \cup L \in
\mathcal{K}^n$, \linebreak the map $\Psi$ is a Minkowski
valuation. Thus, it follows from the equivariance properties of
the Steiner point map $s$ that the map $\Psi - s$ is a
translation invariant, continuous Minkowski valuation that is
equivariant with respect to rotations and homogeneous of degree
one.

Apart from constructing a large class of non-trivial examples, the
main purpose of \cite{Sch74a} was to find reasonable additional
assumptions to single out suitable combinations of dilatations
and reflections among the class of Minkowski endomorphisms. One of
the obtained results was the following, see \cite[Theorem 1.8 (b)]{Sch74a},
which will be used below.

\A

\noindent{\bf Theorem 3.} {\em Let $\Psi:{\cal K}^n \to {\cal
K}^n$ be an endomorphism. If the image under $\Psi$ of some convex
body is a segment, then there are constants $c_2, c_3 \geq 0$
such that}
\[\Psi K =c_2[K-s(K)]+c_3[-K+s(K)]+s(K) \qquad \mbox{for } K \in
\mathcal{K}^n. \]

Motivated by the work of Schneider, the second author \cite{Schu06} recently
investigated continuous operators $\Psi:
\mathcal{K}^n \rightarrow \mathcal{K}^n$, called Blaschke
Minkowski homomorphisms, with the following properties: \\[2mm]
(a) $\Psi$ is Blaschke Minkowski additive, i.e., $\Psi(K\,
\#\,L)=\Psi K + \Psi L$ for $K, L \in \mathcal{K}^n.$\\[2mm]
(b) $\Psi$ is translation invariant and equivariant with respect to rotations.\\[2mm]
Here, the Blaschke sum $K \, \# \, L$ of the convex
bodies $K, L \in \mathcal{K}^n$ is the convex body with $S_{n-1}(K
\, \# \, L,\cdot)=S_{n-1}(K,\cdot)+S_{n-1}(L,\cdot)$ and, say, the
Steiner point at the origin. Property (a) and the continuity of
$\Psi$ imply that $\Psi$ is homogeneous of degree $n-1$.
\linebreak Moreover, since $(K \cup L)\,\#\,(K \cap L)=K\,\#\,L$
whenever $K, L, K \cup L \in \mathcal{K}^n$, the map $\Psi$ is a
Minkowski valuation. A result of McMullen \cite{McM80} on
continuous, translation invariant real valued valuations implies
(compare the proof of Theorem 1.2 in \cite{Schu06}) that Blaschke
Minkowski homomorphisms are precisely the continuous, translation
invariant valuations, homogeneous of degree $n-1$, that are
equivariant with respect to rotations.

Among other results, the following, essentially unique,
representation of Blaschke Minkowski homomorphisms $\Psi$ was
obtained in \cite[Theorem 1.2 and Lemma 4.6]{Schu06}:
\begin{equation}\label{BMhomo}
h(\Psi K,u) =\int_{S^{n-1}}\left[p(\langle u,v \rangle)+ q(\langle
u,v \rangle)\right]\,\D S_{n-1}(K,v), \qquad u \in S^{n-1},
\end{equation}
for $K\in{\cal K}^n$, where $p,q$ are continuous functions on $[-1,1]$, $p$ is even,
$q$ is odd, and $p(\langle\cdot,v\rangle)$ is the restriction of a
support function to $S^{n-1}$. The following result is a version
of Theorem 5.3 in \cite{Schu06}, it corresponds to Theorem 3 for
Blaschke Minkowski homomorphisms.

\A

\noindent{\bf Theorem 4.} {\em Let $\Psi:{\cal K}^n \to {\cal
K}^n$ be a Blaschke Minkowski homomorphism. If the image under
$\Psi$ of some convex body $M$ of dimension at least $n - 1$ is a
polytope, then $\Psi K =c_1\Pi K$ for each centrally symmetric $K
\in \mathcal{K}^n$, where $c_1 \geq 0$ is a real constant.}

\A

We note here that Theorem 5.3 in \cite{Schu06} was formulated for
an $n$-dimensional body $M$ and a certain class of Blaschke
Minkowski homomorphisms, but the proof needs only minor
modifications to give the result stated as Theorem 4; the only
requirement is that the support of $S_{n-1}(M,\cdot)$ is not
empty, which is satisfied if $\dim M = n-1$.

Continuous, translation invariant and rotation equivariant valuations
$\Phi:{\cal K}^n \to {\cal K}^n$ which are homogeneous of some degree
$j\in\{2,\dots,n-2\}$, have not been much investigated. Examples are the
mappings $\Pi_j$ defined by
\[ h(\Pi_j K,u) = \frac{1}{2}\int_{S^{n-1}}|\langle u,v\rangle|\,\D S_j(K,v) \qquad\mbox{for }u\in S^{n-1},\]
where $S_j(K,\cdot)$ is the area measure of order $j$ of the convex body $K$.
The body $\Pi_jK$ is known as the projection body of order $j$ of the convex
body $K$; see \cite[p. 161]{Gar95}. These examples can be generalized considerably.

The proofs of Theorems 1 and 2 will make use of the
following generalization of Theorems 3 and 4, concerning
homogeneous valuations of arbitrary degrees.

\A

\noindent{\bf Theorem 5.} {\em Let $n\ge 3$. Let $\Phi:{\cal K}^n \to {\cal
K}^n$ be a continuous, translation invariant and rotation
equivariant valuation
which is homogeneous of degree $j$ and maps some convex body of dimension $j$ to a polytope.} \\[2mm]
(a) {\em If $j=n$, then $\Phi=\{0\}$.}\\[2mm]
(b) {\em If $j=n-1$, then $\Phi=c_1\Pi$ with a constant $c_1\ge 0$.}\\[2mm]
(c) {\em If $j\in\{2,\dots,n-2\}$, then $\Phi=\{0\}$.}\\[2mm]
(d) {\em If $j=1$, then $\Phi = c_2{\cal I}+ c_3(-{\cal I})$ with constants $c_2,c_3\ge 0$.}\\[2mm]
(e) {\em If $j=0$, then $\Phi=\{0\}$.}

\A

\noindent{\em Proof.} (a) Let $j=n$. For $u\in\R^n$, the
function $\varphi$ defined by $\varphi(K)=h(\Phi K,u)$ satisfies
the assumptions of Lemma 1. It follows that $h(\Phi K,u) =
f(u)V_n(K)$ for $K\in{\cal K}^n$. This defines a
function $f$ on $\R^n$. For $\vartheta\in
SO(n)$ we have
\[ f(\vartheta u) V_n( K) =f(\vartheta
u)V_n(\vartheta K)= h(\Phi \vartheta K,\vartheta u)=h(\vartheta
\Phi K,\vartheta u) =h(\Phi K,u) = f(u)V_n(K),\]
hence $f(u)=a\|u\|$ and thus $\Phi K =a V_n(K)B^n$ with a constant $a$,
where $B^n$ denotes the unit ball. Inserting for $K$ an
$n$-dimensional convex body for which $\Phi K$ is a polytope, we
get that $aB^n$ is a polytope. This is only
possible if $a=0$ and hence $f=0$. This proves part (a).

(b) Let $j=n-1$. Since $\Phi$ is homogeneous of degree $n-1$, it
is a Blaschke Minkowski homomorphism. Since $\Phi$ maps some
convex body $M$ of dimension $n-1$ to a polytope, Theorem 4 yields
that $\Phi K=c\Pi K$ holds for every centrally symmetric convex
body $K$, where $c \geq 0$ is a constant.

Let $K\in{\cal K}^n$, and let  $K \, \, \# \, \, (-K)$ be the Blaschke sum of $K$
and $-K$. Using (\ref{BMhomo}), it is easy to see that $\Phi$
commutes with the reflection in the origin. Thus, it follows from
the fact that $K\,\, \# \, \,(-K)$ is centrally symmetric that
\[\Phi K + (-\Phi K) =\Phi K + \Phi(- K) = \Phi(K\,\, \# \, \,(-K))
=c\Pi(K\,\, \# \, \,(-K)).\]
Suppose, first, that $c=0$. Then $\Phi K$ is one-pointed, say
$\Phi K = \{t(K)\}$. The map $t+s:{\cal K}^n\to\R^n$, where $s$ is the
Steiner point map, is a continuous valuation which is equivariant under
translations and rotations. From the characterization of the Steiner point
mentioned after (\ref{defsteinerpoint}), we obtain that $t=0$, hence
$\Phi K=\{0\} = c\Pi K$ for $K\in{\cal K}^n$.

Let $c>0$. Let $B$ be an $n$-dimensional polytope with the property that any
three of the outer unit normal vectors of its facets are linearly
independent. Writing $\Phi B=:Q$, we get $Q + (-Q) = c\Pi(B\,\, \# \, \,(-B))=:Z$
and, clearly, $c\Pi B = \frac{1}{2}Z$. Then
$Q$ is a summand of the polytope $Z$ and is, therefore, itself a
polytope. Let $F$ be a two-dimensional face of $Q$, and let $u\in
S^{n-1}$ be such that $F=F(Q,u)$, the face of $Q$ with outer
normal vector $u$. Then $F(Q,u)$ is a summand of $F(Z,u)$. Since
the normal cone of $Q$ at $F$ has dimension $n-2$ and hence
cannot be covered by normal cones of $Z$ at faces of dimensions
larger than $2$, we can choose $u$ in such a way that $F(Z,u)$ is
a two-face of $Z$. Due to the way how the directions of the edges
of $Z$ are determined by the facet normals of $B$, the assumption
that any three normal vectors of the facets of $B$ are linearly
independent implies that $F(Z,u)$ is a parallelogram. Every
summand of a parallelogram is either a parallelogram or a segment
or a singleton. In any case, we deduce that $F(Q,u)$ is centrally
symmetric. Since $F$ was an arbitrary two-face of $Q$, the
polytope $Q$ is a zonotope (\cite[Theorem 3.5.1]{Sch93}). In
particular, $Q$ is centrally symmetric. (This holds also if $Q$
has no two-faces.) Therefore, $Q=\frac{1}{2}Z+t(B)$, and thus
$\Phi B=c\Pi B+t(B)$, with some translation vector $t(B)$. Using
(\ref{defsteinerpoint}), (\ref{BMhomo}) and the fact that the
center of mass of $S_{n-1}(B,\cdot)$ is the origin, we obtain
$s(\Phi B)=0=s(c\Pi B)$ and thus $t(B)=0$. Since any
convex body can be approximated by polytopes $B$ satisfying the
assumption on the facet normals, and since $\Phi$ and $\Pi$ are
continuous, we deduce that $\Phi K=c\Pi K$ for all convex bodies
$K\in{\cal K}^n$. This completes the proof of part (b).

(c) Let $j\in\{2,\dots,n-2\}$.  We choose linear subspaces
$E\subset U\subset\R^n$ with $\dim E=j$ and $\dim U=
j+1$. Let $\pi_U:\R^n\to U$ denote the orthogonal
projection, and define a map $\Psi:{\cal K}(U)\to {\cal K}(U)$ by
\[ \Psi K := \pi_U \Phi K \qquad\mbox{for }K\in {\cal K}(U).\]
Then $\Psi$ is a continuous valuation on ${\cal K}(U)$, it is
invariant under the translations of $U$ into itself and
equivariant under the rotations in $SO(U)$. By assumption, $\Phi$
maps some $j$-dimensional convex body to a polytope. By the
translation invariance and rotation equivariance of $\Phi$, there
also exists such a body that is contained in $U$. Now we can
apply Part (b) of Theorem 5, with $\R^n$ replaced by
$U$. It follows that $\Psi = c\Pi^U$, where $\Pi^U$ denotes the
projection body operator in $U$, and $c\ge 0$ is a constant
(depending on $U$). Let $K\subset E$ be any $j$-dimensional
convex body. Then $\Pi^U K=: S$ is a (non-degenerate) segment in $U$, centered at $0$
and orthogonal to $E$, thus
\[ \pi_U \Phi K = cS.\]
In particular, the orthogonal projection of $\Phi K$ to $E$ is
equal to $\{0\}$. Therefore, $\Phi K$ is contained in the orthogonal
complement $E^{\perp}$ (with respect to $\R^n$) of $E$. Every rotation of ${\mathbb
R}^n$ that leaves $E$ pointwise fixed maps $\Phi K$ to itself,
hence $\Phi K$ is a ball with center $0$ and dimension $n-j\ge 2$ or dimension
zero. But for a suitable $j$-dimensional body $M\subset E$, the
set $\Phi M$ is also a polytope. It follows that $\Phi M=\{0\}$.
This implies that $c=0$. In particular, we have $\Phi K =\{0\}$
for every $K \in \mathcal{K}(E)$. Here $E$ can be any
$j$-dimensional subspace, hence $\Phi K=\{0\}$ holds for all
$K\in{\cal K}^n$ with $\dim K \leq j$. An application of Lemma 5 now
completes the proof of part (c).

(d) Let $j=1$. The continuous, translation invariant valuation
$\Phi$ is homogeneous of degree one and hence Minkowski additive
(see \cite{Had57} or \cite{McM77}). The map defined by $K\mapsto
\Phi K+s(K)$ is an endomorphism of ${\cal K}^n$ in the sense of
\cite{Sch74a}. A one-dimensional convex body, that is, a segment,
is a polytope and hence is mapped by $\Phi$ to a polytope. Since
the image has rotational symmetry and $n\ge 3$, it can only be a segment. From
Theorem 3 we can now conclude that
\[ \Phi K = c_2 [K -s(K)]+ c_3[-K+s(K)] \qquad \mbox{for }K\in{\cal K}^n,\]
with constants $c_2,c_3\ge 0$.

(e) Let $j=0$. A continuous translation invariant valuation which
is homogeneous of degree zero is constant, hence for any
$u\in\R^n$ we get $h(\Phi K,u)=f(u)$. As in the proof of
case (a), we obtain $\Phi K^n=aB^n$ with a constant $a\ge 0$.
Choosing for $K$ a one-pointed set for which $\Phi K$ is a
polytope, we get
$a=0$ and hence $\Phi=\{0\}$. This completes the proof of Theorem 5.
\hfill  $\Box$

\A

Note that the proof of Theorem 5 (b) also holds under the
modified assumption that some $n$-dimensional convex body is
mapped to a polytope, which leads to the following generalization
of Theorem 4:

\A

\noindent{\bf Corollary 1.} {\em Let $\Psi:{\cal K}^n \to {\cal
K}^n$ be a Blaschke Minkowski homomorphism. If the image under
$\Psi$ of some convex body $M$ of dimension at least $n - 1$ is a
polytope, then $\Psi K =c \Pi K$, where $c \geq 0$ is a real
constant.}

\section{Proof of Theorem 1}

\noindent We assume that $\Phi:{\cal K}^n\to{\cal K}^n$
satisfies the assumptions of Theorem 1. Let $u
\in\R^n$. By Lemma 2, the real valued valuation $K
\mapsto h(\Phi K,u)$ has a decomposition
\begin{equation}\label{5}
h(\Phi K,u) = \sum_{i=0}^n f_i(K,u), \qquad K\in{\cal K}^n,
\end{equation}
where $f_i(\cdot,u)$ is a continuous translation invariant
valuation that is homogeneous of degree $i$. In (\ref{5}), we
replace $K$ by $mK$ for $m=1,2,\dots,n+1$. The resulting system
of linear equations,
\[ h(\Phi mK,u) = \sum_{i=0}^n m^i f_i(K,u), \qquad m=1,\dots,n+1, \]
can be solved to give representations
\[ f_j(K,u) = \sum_{m=1}^{n+1} a_{jm} h(\Phi mK,u), \qquad j=0,\dots,n, \]
with coefficients $a_{jm}$ depending only on $j$ and $m$. From
this representation we read off the following:\\[2mm]
(a) For each rotation $\vartheta\in SO(n)$ we have $f_j(\vartheta
K,u)=f_j(K,\vartheta^{-1}u)$.\\[2mm]
(b) The function $f_j(K,\cdot)$ is positively homogeneous.\\[2mm]
(c) If $K$ is a polytope, then the function $f_j(K,\cdot)$ is
piecewise linear.

\A

We do not know, at this point, whether each function $f_i(K,\cdot)$ is
a support function; only the following can be shown.

\A

\noindent{\bf Lemma 6.} {\em Suppose that the convex body
$K\in{\cal K}^n$ satisfies
\begin{equation}\label{6}
h(\Phi \lambda K,\cdot) = \sum_{i=k}^l f_i(\lambda K,\cdot)
\end{equation}
for $\lambda>0$, with some $k,l\in\{0,\dots,n\}$, $k\le l$. Then
there exist convex bodies $\Phi_k K, \Phi_l K$ such that
$h(\Phi_k K,\cdot)=f_k(K,\cdot)$ and $h(\Phi_l
K,\cdot)=f_l(K,\cdot)$. If $\Phi \lambda K$ is a polytope for
$\lambda >0$, then $\Phi_k K$ and $\Phi_l K$ are polytopes.

If $E\subseteq \R^n$ is a linear subspace and\/} (\ref{6}) {\em holds for all $K\in{\cal K}(E)$, then the
maps $\Phi_k,\Phi_l:{\cal K}(E)\to{\cal K}^n$ defined in this way
are continuous valuations, invariant under translations and equivariant
under rotations of $E$ into itself, and homo{\-}geneous of degrees
$k$ and $l$, respectively.}

\A

\noindent{\em Proof.} Let $u_1,u_2\in\R^n$ and
$\lambda>0$. Since $h(\Phi\lambda K,\cdot)$ is sublinear,
(\ref{6}) yields
\begin{eqnarray*}
0 &\ge& h(\Phi \lambda K, u_1+u_2) - h(\Phi \lambda K, u_1) - h(\Phi \lambda K, u_2)\\
&=& \sum_{i=k}^l \lambda^i[f_i(K,u_1+u_2)-f_i(K,u_1)-f_i(K,u_2)].
\end{eqnarray*}
Dividing by $\lambda^k$ and letting $\lambda$ tend to zero, we
see that the function $f_k(K,\cdot)$ is sub{\-}linear. Being
positively homogeneous, it is a support function, hence there
exists a convex body $\Phi_kK$ with $f_k(K,\cdot)=h(\Phi_k
K,\cdot)$. If all bodies $\Phi\lambda K$ are polytopes, then
$\Phi_k K$ is a polytope, since $h(\Phi_k K,\cdot)$ is piecewise
linear.

Similarly, dividing by $\lambda^{l}$ and letting $\lambda$ tend
to infinity, we obtain that $f_l(K,\cdot)$ is sublinear and hence
$f_l(K,\cdot)=h(\Phi_l K,\cdot)$ with a convex body $\Phi_lK$.
The remaining assertions are clear.  \hfill  $\Box$

\A

First we apply Lemma 6 with $k=0$ and $l=n$ (and $E=\R^n$). Theorem 5 (a) and (e)
implies that $\Phi_0 =\{0\}$ and $\Phi_n=\{0\}$, hence
\begin{equation}\label{8}
h(\Phi K,\cdot) = \sum_{i=1}^{n-1} f_i(K,\cdot) \qquad \mbox{for
}K\in{\cal K}^n.
\end{equation}
Now Lemma 6 with $l=n-1$ yields the existence of a map
$\Phi_{n-1}:{\cal K}^n\to {\cal K}^n$ which is a continuous,
translation invariant and rotation equivariant valuation,
homogeneous of degree $n-1$, and satisfying
$h(\Phi_{n-1}K,\cdot)=f_{n-1}(K,\cdot)$. If $K$ is a polytope,
then $\Phi_{n-1}K$ is a polytope. Theorem 5 (b) shows that
$\Phi_{n-1} = c_1\Pi$ with a constant $c_1\ge 0$.

Similarly, we conclude from (\ref{8}) that
$f_1(K,\cdot)=h(\Phi_1K,\cdot)$ with a continuous, translation
invariant and rotation equivariant valuation $\Phi_{1}:{\cal
K}^n\to {\cal K}^n$ which is homogeneous of degree one, and that
$\Phi_1 K$ is a polytope if $K$ is a segment. From Theorem 5 (d)
we obtain that $ \Phi_1 = c_2 {\cal I}+ c_3(-{\cal I})$
with constants $c_2,c_3\ge 0$. Therefore, (\ref{8}) can be replaced by
\begin{equation}\label{9}
h(\Phi K,\cdot) = c_1 h(\Pi K,\cdot)+\sum_{i=2}^{n-2} f_i(K,\cdot) + c_2 h({\cal I}K ,\cdot)+ c_3h(-{\cal I}K,\cdot).
\end{equation}
This finishes the proof if $n=3$. We assume, therefore, that
$n\ge 4$. We have to show that the remaining functions
$f_i(K,\cdot)$ are zero.

Let $j\in\{2,\dots,n-2\}$. We choose a $j$-dimensional linear
subspace $E\subset\R^n$. Let $K\in{\cal K}(E)$. Since a
continuous, translation invariant valuation that is homogeneous of
degree $i$ vanishes on convex bodies of dimension smaller than
$i$, we have
\[ h(\Phi \lambda K,\cdot) = \sum_{i=1}^j f_i(\lambda K,\cdot)\]
for $\lambda>0$. By Lemma 6, there is a convex body $\Phi_j K$
with $h(\Phi_j K,\cdot) = f_j(K,\cdot)$, and if $K$ is a
polytope, then $\Phi_j K$ is a polytope. For $u\in\R^n$,
Lemma 1 gives $h(\Phi_j K,u) = f(u) V_j(K)$, with a
function $f$ on $\R^n$. Taking
for $K$ a $j$-dimensional polytope in $E$, we see that $f$ is
the support function of a polytope,
say $P$. By the rotation equivariance, $P$ is invariant under the
rotations mapping $E$ into itself and keeping $E^{\perp}$
pointwise fixed. Therefore, the projection $\pi_E P$ is a ball in
$E$, centered at $0$. It can only have radius zero, hence $P$ is contained in
$E^{\perp}$. Every rotation of $\R^n$ that
leaves $E$ pointwise fixed maps $P$ to itself, hence $P$ is a centered ball
of dimension $n-j\ge 2$ or of dimension zero. We deduce that
$P=\{0\}$.

We have shown that
\begin{equation}\label{77}
f_j(K,\cdot)=0 \mbox{ whenever } \dim K= j,\; j=2,\dots,n-2.
\end{equation}
From Lemma 4 (a) we conclude that
\begin{equation}\label{10}
f_j(K,\cdot) +f_j(-K,\cdot)=0 \qquad \mbox{ for all }K\in{\cal K}^n,\; j=2,\dots,n-2.
\end{equation}

Now let $E\subset\R^n$ be an $(n-1)$-dimensional linear subspace.
Define $\Psi: \mathcal{K}(E) \rightarrow
\mathcal{K}(E)$ by $\Psi K=\pi_E \Phi K$ for $K \in
\mathcal{K}(E)$. Then $\Psi$ is a continuous valuation,
invariant under the translations of $E$ into itself and
equivariant under the rotations of $SO(E)$. It maps polytopes to
polytopes. Let $K\in{\cal K}^n$. The support function, on $E$, of $\Psi K$ is the
restriction of $h(\Phi K,\cdot)$ to $E$. Hence, it follows from
(\ref{8}) that
\[h(\Psi K,u)=\sum \limits_{i=1}^{n-1} f_i(K,u) \qquad \mbox{for }u \in E. \]
On the other hand, from the result (\ref{9}), applied in $E$, we
have
\[h(\Psi K,u)=c_Eh(\Pi^E K,u)+\sum\limits_{i=1}^{n-3}g_i(K,u) \qquad \mbox{for }u \in E, \]
where $\Pi^E$ is the projection body operator in $E$ and
$g_i(\cdot,u)$ is homogeneous of degree $i$. By homogeneity, we
must have
\[c_Eh(\Pi^EK,u)=f_{n-2}(K,u) \qquad \mbox{for }u \in E. \]
Let $K \subset E$ be an $(n-1)$-dimensional centrally symmetric
body. Since $\Pi^EK \neq \{0\}$, we deduce from (\ref{10}) that
$c_E=0$. This yields
\begin{equation} \label{15}
f_{n-2}(K,u)=0 \qquad \mbox{for all }K \in \mathcal{K}(E),\; u \in E.
\end{equation}

Let $e$ be one of the unit normal vectors of $E$, and let
$K\in{\cal K}(E)$. Then (\ref{9}) can be written in the form
\begin{equation}\label{11}
h(\Phi K,\cdot) = c_1 V_{n-1}(K)|\langle e,\cdot\rangle|+\sum_{i=1}^{n-2} f_i(K,\cdot).
\end{equation}

Let $H^+_e:=\{u\in\R^n: \langle e,u\rangle \ge 0\}$ and
$H^-_e=-H^+_e$. For $u_1,u_2\in H^+_e$ and $\lambda>0$, we have
\begin{eqnarray*}
0 &\ge& h(\Phi \lambda K, u_1+u_2) - h(\Phi \lambda K, u_1) - h(\Phi \lambda K, u_2)\\
&=& \sum_{i=1}^{n-2} \lambda^i[f_i(K,u_1+u_2)-f_i(K,u_1)-f_i(K,u_2)],
\end{eqnarray*}
from which we obtain
\[ f_{n-2}(K,u_1+u_2) \le f_{n-2}(K,u_1) + f_{n-2}(K,u_2).\]
We replace $K$ by $-K$ and use (\ref{10}). Together with the
preceding inequality, this yields
\[ f_{n-2}(K,u_1+u_2) = f_{n-2}(K,u_1) + f_{n-2}(K,u_2).\]
Since this holds for all $u_1,u_2\in H^+_e$ and since
$f_{n-2}(K,\cdot)$ is positively homogeneous, we conclude that
$f_{n-2}(K,\cdot)$ is linear on $H^+_e$, thus there is a vector
$x_K$ such that
\[ f_{n-2}(K,u) = \langle x_K,u \rangle \qquad\mbox{for }u\in H^+_e.\]
Similarly, there is a vector $y_K$ such that
\[ f_{n-2}(K,u) = \langle y_K,u \rangle \qquad\mbox{for }u\in H^-_e.\]
By (\ref{15}), for any vector $v\perp e$ we have $\langle x_K,v
\rangle = f_{n-2}(K,v)= 0$, and analogously $\langle y_K,v \rangle =0$, hence $x_K$
and $y_K$ are parallel to $e$. Thus,
\[f_{n-2}(K,u)=\langle e,u \rangle \varphi_{n-2}^{\pm}(K) \qquad \mbox{for }u \in H_e^{\pm}, \]
where $\varphi_{n-2}^+,\varphi_{n-2}^-: \mathcal{K}(E) \rightarrow \mathbb{R}$ are
translation invariant, continuous valuations,
homogeneous of degree $n-2$, and by (\ref{77}) they are simple.
Moreover, for every rotation $\vartheta \in SO(E)$ fixing $e$ we have
\[\varphi_{n-2}^{+}(\vartheta K)=f(\vartheta K,e)=f(K,e)= \varphi_{n-2}^{+}( K). \]
Thus, $\varphi_{n-2}^+$ is also invariant with respect to rotations
of $E$ into itself. By Lemma 3, applied in $E$, this is only possible if
$\varphi_{n-2}^+=0$. Similarly, we obtain $\varphi_{n-2}^-=0$ and thus
$f_{n-2}(K,\cdot)=0$ for all $K\in \mathcal{K}(E)$. Since $E$ was
an arbitrary $(n-1)$-dimensional subspace, we have
$f_{n-2}(K,\cdot)=0$ for arbitrary convex bodies $K$ with $\dim
K\le n-1$. Now Lemma 4 (b) yields $f_{n-2}(K,\cdot)=0$ for all
convex bodies $K \in \mathcal{K}^n$.

Let $j\in\{2,\dots,n-3\}$ and suppose it has already been proved that
\[ f_i(K,\cdot) =0 \qquad\mbox{for all }K\in{\cal K}^n,\; i=j+1,\dots,n-2. \]
We choose a $(j+1)$-dimensional linear
subspace $E\subset\R^n$. For all $K\in{\cal K}(E)$ we have
$\Pi K=\{0\}$ and hence, by (\ref{9}),
\[ h(\Phi K,\cdot) = \sum_{i=1}^j f_i(K,\cdot).\]
By Lemma 6, the valuation $\Phi_j:{\cal K}(E)\to{\cal K}^n$ with
$h(\Phi_j K,\cdot) = f_j(K,\cdot)$ is defined. By (\ref{77}),
it satisfies $\Phi_j K=0$ if $\dim K=j$. The proof of Lemma 5 shows that
$\Phi_j =\{0\}$. Thus, $f_j(K,\cdot)=0$ whenever
$\dim K \le j+1$. Now Lemma 4 (b) yields
\[ f_j(K,\cdot) =0 \qquad\mbox{for all }K\in{\cal K}^n. \]
In this way we continue, until we obtain $f_2(K,\cdot)=0$ for all $K\in{\cal K}^n$.
This completes the proof of Theorem 1. \hfill  $\Box$

\section{Proof of Theorem 2}

We assume that the assumptions of Theorem 2 are
satisfied. As in the proof of Theorem 1, relation (\ref{5}) holds
(but $f_j(K,\cdot)$ need not have property (c) listed there), thus
\[ h(\Phi \lambda K,\cdot) = \sum_{j=0}^n \lambda^j f_j(K,\cdot)\qquad\mbox{for }K\in{\cal K}^n, \]
for all $\lambda>0$. If $\dim K\le n-2$, then $\Phi \lambda
K=\{0\}$ by assumption and continuity, and we deduce that
$f_j(K,\cdot) =0$ for $j=0,\dots,n$.

We assert that $f_j(K,\cdot)=0$ for $K\in{\cal K}^n$ and for
$j=0,\dots,n-2$. For $j=0$, this follows from $f_0(K,u) = a\|u\|$
(obtained as before), by inserting a convex body $K$ for which
$f_0(K,\cdot)=0$. Suppose that $j\in\{1,\dots,n-2\}$ and that
$f_i(K,\cdot)=0$ for all $K\in{\cal K}^n$ has been proved for $i<j$. Then
\[ h(\Phi K,\cdot) = \sum_{i=j}^n f_i(K,\cdot) \qquad \mbox{for }K\in{\cal K}^n.\]
By Lemma 6, the map $\Phi_j:{\cal K}^n\to{\cal K}^n$ is defined. It maps bodies of
dimension smaller than $n-1$ to $\{0\}$; in particular, $\Phi_jK=0$  if $\dim K=j$.
Lemma 5 shows that $\Phi_j=\{0\}$. Thus, (\ref{5}) reduces to
\[h(\Phi K,\cdot) = f_n(K,\cdot) + f_{n-1}(K,\cdot)= h(\Phi_n
K,\cdot) + h(\Phi_{n-1} K,\cdot)\] for $K\in{\cal K}^n$, where we
have already inserted the homogeneous valuations
$\Phi_n,\Phi_{n-1}$ that exist by Lemma 6.

By Lemma 1, $h(\Phi_n K,u) = g(u)V_n(K)$ with some function $g$.
By the rotation equivariance of $\Phi$, we obtain $\Phi_nK =
aV_n(K)B^n$, with a constant $a\ge 0$. By assumption, there exists
an $n$-dimensional convex body $M$ such that $\Phi M$ is a
polytope $P$. This gives $P=aV_n(M) B^n + \Phi_{n-1}M$, where
$B^n$ denotes the unit ball. Since the polytope $P$ cannot have a
ball with positive radius as a summand, this yields $a=0$.
Therefore, $\Phi = \Phi_{n-1}$. Now Corollary 1 completes the
proof of Theorem 2. \hfill  $\Box$

\section{Open Problems}

A positive answer to the following problem would
simplify the proof of Theorem 1 considerably:

\A

\noindent{\bf Problem 1.} {\em Let $n\ge 3$. Let $\Phi:{\cal K}^n
\to {\cal K}^n$ be a continuous, translation invariant and
rotation equivariant valuation. Is there a $($unique$)$
representation of $\Phi$ of the form
\[ \Phi = \Phi_0 + \dots + \Phi_n,\]
where $\Phi_j: {\cal K}^n \to \mathcal{K}^n$ is a continuous,
translation invariant and rotation equivariant valuation which is
homogeneous of degree $j$?}

\A

As remarked in the introduction, the assumptions of Theorem 1 are
stronger than necessary. Theorems 2 and 5 lead to the following
question.

\A

\noindent{\bf Problem 2.} {\em Let $n\ge 3$. Let $\Phi:{\cal K}^n
\to {\cal K}^n$ be a continuous, translation invariant and
rotation equivariant valuation. Assume that $\Phi$ maps some
$n$-dimensional convex body to a polytope. Are there constants
$c_1,c_2,c_3\ge 0$ such that}
\[ \Phi = c_1\Pi + c_2 {\cal I}+c_3(-{\cal I})? \]

\A

\noindent{\bf \Large Acknowledgments} \\[5mm]
The authors would like to thank Monika Ludwig for helpful
remarks. This work was supported by the European Network PHD,
FP6 Marie Curie Actions, RTN, Contract MCRN-2004-511953. The
second named author was also supported by the Austrian Science Fund
(FWF), within the project ``Affinely associated bodies'', Project
Number: P16547-N12.

\vspace{5mm}

\noindent Rolf Schneider: Mathematisches Institut,Universit\"at Freiburg, Eckerstra{\ss}e 1, \\79104 Freiburg i. Br., Germany\\
E-mail address: rolf.schneider@math.uni-freiburg.de\\[5mm]
Franz E. Schuster: Institut f\"ur Diskrete Mathematik und Geometrie,\\Technische Universit\"at Wien, Wiedner Hauptstra{\ss}e 8-10/104, 1040 Wien, Austria\\
E-mail address: franz.schuster@tuwien.ac.at


\begin{thebibliography}{10}


\bibitem{Ales99} S. Alesker, {\em Continuous rotation invariant
valuations on convex sets}, Ann. of Math (2) {\bf 149} (1999),
977--1005.

\bibitem{Ales01} S. Alesker, {\em Description of translation
invariant valuations on convex sets with solution of P.
McMullen's conjecture}, Geom. Funct. Anal. {\bf 11} (2001),
244--272.

\bibitem{bolk69} E. D. Bolker, {\em A class of convex bodies},
Trans. Amer. Math. Soc. {\bf 145} (1969), 323--345.

\bibitem{bourg88} J. Bourgain and J. Lindenstrauss, {\em Projection bodies}, in: Geometric Aspects of
Functional Analysis (1986/87), Lecture Notes in Math. 1317, Springer, Berlin, 1988, pp. 250--270.

\bibitem{Gar95} R. J. Gardner, {\em Geometric Tomography}, Cambridge University Press, Cambridge, 1995.

\bibitem{GW93} P. Goodey and W. Weil, {\em Zonoids and generalisations}, in: P. M. Gruber, J. M. Wills (Eds.), Handbook of Convex Geometry, Vol. B,
North-Holland, Amsterdam, 1993, pp. 1297--1326.

\bibitem{grin99} E. Grinberg and G. Zhang,
{\em Convolutions, transforms, and convex bodies}, Proc. London Math. Soc. (3) {\bf 78} (1999), 77--115.

\bibitem{Had57} H. Hadwiger, {\em Vorlesungen \"uber Inhalt, Oberfl\"ache und Isoperimetrie}, Springer, Berlin, 1957.

\bibitem{Kid05} M. Kiderlen, {\em Blaschke- and
Minkowski-endomorphisms of convex bodies}, Trans. Amer. Math. Soc., to appear.

\bibitem{Kla95} D. A. Klain, {\em A short proof of Hadwiger's characterization theorem}, Mathematika {\bf 42} (1995), 329--339.

\bibitem{Kla99} D. A. Klain, {\em Even valuations on convex bodies}, Trans. Amer. Math. Soc. {\bf 352} (1999), 71--93.

\bibitem{KlaRo97} D. A. Klain and G. Rota, {\em Introduction to
Geometric Probability}, Cambridge University Press, Cambridge, 1997.

\bibitem{Lud02} M. Ludwig, {\em Projection bodies and valuations}, Adv. Math. {\bf 172} (2002), 158--168.

\bibitem{Lud03} M. Ludwig, {\em Ellipsoids and matrix valued
valuations}, Duke Math. J. {\bf 119} (2003), 159--188.

\bibitem{Lud05} M. Ludwig, {\em Minkowski valuations}, Trans. Amer. Math. Soc. {\bf 357} (2005), 4191--4213.

\bibitem{Lud06} M. Ludwig, {\em Intersection bodies and
valuations}, Amer. J. Math., to appear.

\bibitem{Lud99} M. Ludwig and M. Reitzner, {\em A characterization of
affine surface area}, Adv. Math. {\bf 147} (1999), 138--172.

\bibitem{McM77} P. McMullen, {\em Valuations and Euler-type relations on certain classes of convex polytopes}, Proc. London Math. Soc. {\bf 35} (1977), 113--135.

\bibitem{McM80} P. McMullen, {\em Continuous translation
invariant valuations on the space of compact convex sets}, Arch. Math. {\bf 34} (1980), 377--384.

\bibitem{McM93} P. McMullen, {\em Valuations and dissections}, in: P. M. Gruber, J. M. Wills (Eds.), Handbook of Convex Geometry, Vol. B,
North-Holland, Amsterdam, 1993, pp. 933--990.

\bibitem{MS83} P. McMullen and R. Schneider, {\em Valuations on convex bodies}, in: P. M. Gruber, J. M. Wills. (Eds.), Convexity and Its Applications, Birkh\"auser, Basel, 1983, pp. 170--247.

\bibitem{Min11} H. Minkowski, {\em Theorie der konvexen K\"orper, insbesondere Begr\"undung ihres Ober\-fl\"achen\-begriffs}, Gesammelte Abhandlungen, Vol. II, Teubner, Leipzig, 1911, pp. 131--229.

\bibitem{Sch72} R. Schneider, {\em Kr\"ummungsschwerpunkte konvexer K\"orper, II}, Abh. Math. Sem. Univ. Hamburg {\bf 37} (1972), 204--217.

\bibitem{Sch74a} R. Schneider, {\em Equivariant endomorphisms of the space of convex bodies}, Trans. Amer. Math. Soc. {\bf 194} (1974), 53--78.

\bibitem{Sch74} R. Schneider, {\em Bewegungs\"aquivariante, additive und stetige Transformationen konvexer Bereiche}, Arch. Math. {\bf 25} (1974), 303--312.

\bibitem{Sch93} R. Schneider, {\em Convex Bodies: The Brunn-Minkowski Theory}, Cambridge University Press, Cambridge, 1993.

\bibitem{Sch96} R. Schneider, {\em Simple valuations on convex bodies}, Mathematika {\bf 43} (1996), 32--39.

\bibitem{SW83} R. Schneider and W. Weil, {\em Zonoids and related topics}, in: P. M. Gruber, J. M. Wills. (Eds.), Convexity and Its Applications, Birkh\"auser, Basel, 1983, pp. 296--317.

\bibitem{Schu06} F. E. Schuster, {\em Convolutions and multiplier transformations}, Trans. Amer. Math. Soc., to appear.

\bibitem{Schu06a} F. E. Schuster, {\em Volume inequalities and additive maps of convex bodies}, Mathematika, to appear.

\bibitem{Thom96} A. C. Thompson, {\em Minkowski Geometry}, Cambridge University Press, Cambridge, 1996.

\bibitem{Zieg95} G. M. Ziegler, {\em Lectures on Polytopes}, Springer, New York, 1995.

\end{thebibliography}
\end{document}